\documentclass[10pt]{article}
\usepackage[T2A]{fontenc}
\usepackage[russian,english]{babel}

\usepackage{amssymb,amsmath,amsthm,amsfonts}
\usepackage{bm,bbm}
\usepackage{mathtools}
\usepackage{xcolor}
\usepackage{comment}
\usepackage{todonotes}
\usepackage[top=2cm,bottom=2cm,left=3cm,right=3cm]{geometry}
\usepackage{hyperref}
\hypersetup{
    colorlinks = true,
    citecolor = blue
}
\usepackage{algorithm}
\usepackage{algorithmic}

\newcommand*{\BREAK}{\textbf{break}}

\usepackage{graphicx}
\usepackage{subcaption}

\usepackage[authoryear]{natbib}
\bibpunct{(}{)}{;}{a}{,}{,}
\DeclareMathOperator*{\argmax}{arg\,max}
\DeclareMathOperator*{\argmin}{arg\,min}

\DeclareMathOperator{\dom}{dom}
\DeclarePairedDelimiter{\norm}{\lVert}{\rVert}

\DeclareMathOperator{\sgr}{\partial}

\newcommand*{\tbar}{\bar{t}}
\newcommand*{\fbar}{\bar{f}}

\newcommand*{\R}{\mathbb{R}}
\newcommand*{\N}{\mathbb{N}}
\newcommand*{\eps}{\varepsilon}

\newtheorem{theorem}{Theorem}[section]
\newtheorem{lemma}[theorem]{Lemma}

\begin{document}
\title{Finding Equilibria in the Traffic Assignment Problem with Primal-Dual Gradient Methods for Stable Dynamics Model and Beckmann Model}

\author{Meruza Kubentayeva\\
\normalsize{Institute for Information Transmission Problems RAS, \href{mailto:kubentay@gmail.com}{kubentay@gmail.com}}
\and Alexander Gasnikov\\ 
\normalsize {Moscow Institute of Physics and Technology,} \\
\normalsize {Institute for Information Transmission Problems RAS,} \\
\normalsize {Higher School of Economics, 
\href{mailto:gasnikov@yandex.ru}{gasnikov@yandex.ru}}
}


\date{}
\maketitle

\begin{abstract}
In this paper we consider the application of several gradient methods to the traffic assignment problem: we search equilibria in the stable dynamics model \citep{nesterov2003stationary} and the Beckmann model. Unlike the celebrated Frank--Wolfe algorithm widely used for the Beckmann model, these gradients methods solve the dual problem and then reconstruct a solution to the primal one. We deal with the universal gradient method, the universal method of similar triangles, and the method of weighted dual averages, and estimate their complexity for the problem. Due to the primal-dual nature of these methods, we use a duality gap in a stopping criterion. In particular, we present a novel way to reconstruct admissible flows in the stable dynamics model, which provides us with a computable duality gap.


\end{abstract}

\section{Introduction}

The Beckmann model for searching static traffic equilibria in road networks is among the most widely used models by transportation planners \citep{beckmann1956studies, patriksson2015traffic}. The equilibria found are practical for evaluating the network efficiency and distribution of business centers and residential areas, and establishing urban development plans, etc.
This model introduces a cost function on every link of a transportation network, which defines a dependence of the travel cost on the flow along the link. In practice the BPR functions are usually employed \citep{bpr_functions}:
\begin{equation}
\label{eq:Beckmann_cost_func} 
    \tau_e(f_e) = \tbar_e \left(1 + \rho \left( \frac{f_e}{\fbar_e}\right)^{\frac{1}{\mu}} \right), 
\end{equation}
where $\tbar_e$ are free flow times, and $\fbar_e$ are road capacities of a given network's link $e$. We take these functions with parameters $\rho = 0.15$ and $\mu = 0.25$.

\begin{figure}[H]
    \centering
    \includegraphics[scale = 0.25, trim={0cm 0 0cm 0},clip]{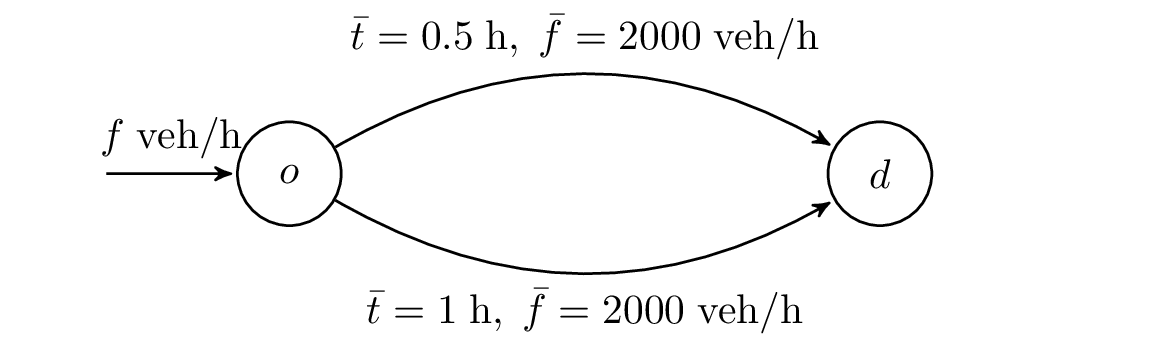}
    \caption{Parallel routes}
    \label{fig:example_pic}
\end{figure}

\cite{nesterov2003stationary} proposed an alternative model called the stable dynamics model, which takes an intermediate place between static and dynamic network assignment models. Namely, its equilibrium can be interpreted as the stationary regime of some dynamic process. Its key assumption is that we no longer introduce a complex dependence of the travel cost on the flow (as in the standard static models), but only pose capacity constraints, i.e.\ the flow value on each link imposes the feasible set of travel times
\begin{equation}
\label{eq:SD_cost_func}
    \tau_e(f_e) = 
    \begin{cases} \tbar_e, & 0 \le f_e < \fbar_e,\\
        \left[ \tbar_e, \infty \right], & f_e = \fbar_e,\\
        +\infty, & f_e > \fbar_e.
    \end{cases}
\end{equation}
Unlike in the Beckmann model, there is no one-to-one correspondence between equilibrium travel times and flows on the links of the network. We can illustrate the difference on a simple example of two parallel routes (Figure~\ref{fig:example_pic}).
Let the input flow take values 1000, 2000, and 3000~veh/h. For the stable dynamics model, in the first and second cases, all drivers choose the upper route; the equilibrium travel time simply equals the upper route's free flow time (0.5~h) in the first case and varies from 0.5~h to 1~h (according to the model) in the second. In the third case, the input flow exceeds the upper route's capacity, so the upper route's flow is 2000~veh/h, the lower one's is 1000~veh/h, and the equilibrium travel time is 1~h. 
All these equilibria can be interpreted as stationary regimes of some dynamic processes, e.g.\ the last case can be viewed as the result of the queue at the beginning of the upper route (since this route's capacity is smaller than the input flow) created by drivers who wanted to take this route until the waiting time plus the route's travel time reached the lower route's travel time \citep{nesterov2003stationary}.
In the Beckmann model equilibria are as follows: \textit{for all three cases} only the upper route is used, and the equilibrium travel times are approximately 0.5~h, 0.6~h, and 0.9~h, respectively.
\cite{chudak2007static} conducted a detailed comparison~--- for large and small networks~--- of equilibria in these two models.

In the Beckmann model, searching equilibria reduces to minimization of a potential function. One of the most popular and effective approaches to solve this problem numerically is the famous Frank--Wolfe method \citep{frank1956algorithm, jaggi2013revisiting} as well as its numerous modifications \citep{fukushima1984modified, leblanc1985improved, arezki1990full, chen2002faster}.

In the case of the stable dynamics model, one cannot directly apply the Frank--Wolfe method. However, an equilibrium can be found as a solution of a pair of primal and dual optimization problems. The same holds also for the Beckmann model, so in both cases we can apply primal-dual (sub)gradient methods.

In this work, we compare several primal-dual gradient methods for searching equilibria in both the Beckmann and the stable dynamics models, namely, the universal gradient method (UGM) \citep{nesterov2015universal}, the universal method of similar triangles (UMST) \citep{gasnikov2018universal}, and the method of weighted dual averages (WDA) \citep{nesterov2009primal}. The main advantage of the above universal methods is an automatic adjustment to a local (H\"older) smoothness of a minimized function, which is especially important since the dual problems we are dealing with are essentially non-smooth.
Due to the primal-dual nature of these methods, one can use an adaptive stopping criterion guaranteeing required accuracy.

The \textbf{main contributions} of this paper include the following:
\begin{itemize}
    \item We propose a novel way to reconstruct admissible flows (i.e.,\ meeting the capacity constraints and induced by flows on the paths) in the stable dynamics model and a novel computable duality gap, which can be used in a stopping criterion.
    \item We provide theoretical upper bounds on the complexity of searching equilibria by the considered algorithms: UMST, UGD, and~WDA.
    \item We conducted numerical experiments comparing these algorithms on the Anaheim transportation network---the source code is available for use and can be found in~\cite{meruzakub}.
\end{itemize}																																		 

The paper is organized as follows.
In Section~\ref{sec:problem_statement} we give a problem statement, define equilibria in the Beckmann and the stable dynamics models and corresponding optimization problems.
Section~\ref{sec:methods} is devoted to the complexity analysis of UGM, UMST, and WDA. We show that the number of iterations required to obtain an $\eps$-solution of primal and dual problems is $O(1 / \eps^2)$ for UGM and  UMST.
In Section~\ref{sec:experiments} results of experiments on Anaheim transportation network are presented.
Finally, some conclusions are drawn in Section~\ref{sec:conclusion}.

\section{Problem statement}
\label{sec:problem_statement}

Let the urban road network be represented by a directed graph $G = ( V, E )$, where vertices $V$ correspond to intersections or centroids \citep{sheffi1985urban} and edges $E$ correspond to roads, respectively.
Suppose we are given the travel demands: namely, let $d_w$(veh/h) be a trip rate for an origin-destination pair $w$ from the set $OD \subseteq \{ w = (i, j) : i \in O, \; j \in D \}$. Here $O \subseteq V$ is the set of all possible origins of trips, and $D \subseteq V$ is the set of destination nodes.
For OD pair $w = (i, j)$ denote by $P_w$ the set of all simple paths from $i$ to $j$. Respectively, $P = \bigcup_{w \in OD} P_w$ is the set of all possible routes for all OD pairs. 
Agents traveling from node $i$ to node $j$ are distributed among paths from $P_w$, i.e.\ for any $p \in P_w$ there is a flow $x_p \in \R_+$ along the path $p$, and $\sum_{p \in P_w} x_p = d_w$.
Flows from vertices from the set $O$ to vertices from the set $D$ create the traffic in the entire network $G$, which can be represented by an element of
\[
X = \Bigl\{x \in \R_{+}^{|P|} : \; \sum_{p \in P_w} x_p = d_w, \; w \in OD \Bigr\}.
\]
Note that the dimension of $X$ can be extremely large: e.g.\ for $n \times n$ Manhattan network $\log |P| = \Omega(n)$.
To describe a state of the network we do not need to know an entire vector $x$, but only flows on arcs:
\[
f_e(x) = \sum_{p \in P} \delta_{e p} x_p \quad \text{for} \quad e \in E,
\]
where $\delta_{e p} = \mathbbm{1}\{e \in p\}$. Let us introduce a matrix $\Theta$ such that $\Theta_{e, p} = \delta_{e p}$ for $e \in E$, $p \in P$, so in vector notation we have $f = \Theta x$. To describe an equilibrium we use both path- and link-based notations $(x, t)$ or $(f, t)$. 

\textbf{Beckmann model.} One of the key ideas behind the Beckmann model is that the cost (e.g.\ travel time, gas expenses, etc.) of passing a link $e$ is the same for all agents and depends solely on the flow $f_e$ along it. In what follows, we denote this cost for a given flow $f_e$ by $t_e = \tau_e(f_e)$.
Another essential point is a behavioral assumption on agents called the first Wardrop's principle: we suppose that each of them knows the state of the whole network and chooses a path $p$ minimizing the total cost
\[
T_p(t) = \sum_{e \in p} t_e.
\]

The cost functions are supposed to be continuous, non-decreasing, and non-negative. Then $(x^*, t^*)$, where $t^* = (t_e^*)_{e \in E}$, is an equilibrium state, i.e.\ it satisfies conditions
\begin{gather*}
    t_e^* = \tau_e(f_e^*), \quad\text{where}\quad f^* = \Theta x^*, \\ 
    x^*_{p_w} > 0 \Longrightarrow T_{p_w}(t^*) = T_w(t^*) = \min_{p \in P_w} T_p(t^*),
\end{gather*}
if and only if $x^*$ is a minimum of the potential function:
\begin{align*}
    \Psi(x) = \sum_{e \in E} \underbrace{\int_{0}^{f_e} \tau_e (z) d z}_{\sigma_e(f_e)} 
    \longrightarrow \min_{f = \Theta x, \; x \in X} \\
    \Longleftrightarrow \Psi(f) = \sum_{e \in E} \sigma_e (f_e) 
    \longrightarrow \min_{f = \Theta x : \; x \in X}, \tag{B}\label{PrimalBeckmann}
\end{align*}
and $t_e^* = \tau_e(f_e^*)$ \citep{beckmann1956studies}.

Another way to find an equilibrium numerically is by solving a dual problem.
According to Theorem~4 from \cite{nesterov2003stationary}, we can construct it in the following way:
\begin{align*}
    \min_{f = \Theta x : \; x \in X} \Psi(f) &= \min_{x\in X, \; f} \left [ \Psi(f) + \sup_{t \in \R^{|E|}} \langle t, \Theta x - f \rangle \right ] = \sup_{t \in \R^{|E|}} \min_{x \in X, \; f} \left [ \Psi(f) + \langle t, \Theta x - f \rangle \right ] \\
    &= \sup_{t \in \R^{|E|}} \left [ - \sum_{e \in E} \max_{f_e} \{t_e f_e - \sigma_e(f_e) \} + \min_{x \in X} \sum_p \sum_{e \in E} t_e \delta_{ep} x_p \right ] \\
    &= \max_{t \in \dom \sigma^*} - \left [ \sum_{e \in E} \sigma_e^*(t_e) - \sum_{w \in OD} d_w T_w(t) \right ] = - \min_{t \ge \bar{t}} Q(t),
\end{align*}
where 
\[
\sigma_e^*(t_e) = \sup_{f_e \ge 0} \{t_e f_e - \sigma_e(f_e) \} 
= \fbar_e \left( \frac{t_e - \tbar_e}{\tbar_e \rho} \right)^{\mu} \frac{\left(t_e - \bar{t}_e \right)}{1 + \mu}
\] 
is the conjugate function of $\sigma_e(f_e)$, $e \in E$.
Finally, we obtain the dual problem, the solution of which is $t^*$:
\begin{equation}\label{DualBeckmann}
    Q(t) = - \sum_{w \in OD} d_w T_w(t) + \sum_{e \in E} \sigma_e^*(t_e) \longrightarrow \min_{t \ge \bar{t}}. \tag{DualB}
\end{equation}

When we search for the solution to this problem  numerically, on every step of an applied method we can reconstruct primal variable $f$ from the current dual variable $t$: $f \in \partial \sum_{w \in OD} d_w T_w(t)$ (see Subsection~\ref{subsec:subgrad}). Then we can use the duality gap~--- which is always nonnegative~--- for the estimation of the method's accuracy: 
\[
\Delta(f, t) = \Psi(f) + Q(t).
\]
It vanishes only at the equilibrium $(f^*, t^*)$. 



\textbf{Stable dynamics model \citep{nesterov2003stationary}.}
An equilibrium state $(x^*, t^*)$ of the stable dynamics model satisfies the next conditions:
\begin{gather*}
    t_e^* \in \tau_e(f_e^*), \\ 
    x^*_{p_w} > 0 \Longrightarrow T_{p_w}(t^*) = T_w(t^*),
\end{gather*}
where $\tau(f)$ is defined earlier by \eqref{eq:SD_cost_func}.
The above formula can be reformulated in terms of an optimization problem:
\begin{align*}
    x^* &= \argmin\limits_{x \in X} \sum_{w \in OD} \sum_{p \in P_w} x_{p} T_p(t^*) \\
    &= \argmin\limits_{x \in X} \sum_{e \in E} t^*_e f_e(x) \\
    &= \argmin_{x \in X} \sum_{e \in E} [t^*_e f_e(x) - \underline{ (t_e^* - \tbar_e) \fbar_e}],
\end{align*}
\begin{align*}
    t_e^* \in \tau_e (f_e^*) 
    \Longleftrightarrow t_e^* &= \argmax_{t_e \ge \tbar_e} t_e (f_e^* - \fbar_e)\\
    &= \argmax_{t_e \ge \tbar_e} [t_e (f_e^* - \fbar_e) + \underline{\tbar_e \fbar_e}].
\end{align*}
Here, we add underlined constant terms to show that the pair $(f^*, t^*)$ is an equilibrium if and only if it is a solution of the saddle-point problem
\begin{equation*}\label{SaddleSD}
\sum_{e \in E} [t_e f_e - (t_e - \tbar_e) \fbar_e] 
\longrightarrow \min_{\substack{f = \Theta x: \\ x \in X}} \max_{t_e \ge \tbar_e}, \tag{SaddleSD}
\end{equation*}
where its primal problem is
\begin{align*}
    &\Psi(x) = \sup_{t_e \ge \tbar_e} \sum_{e \in E} [t_e f_e - (t_e - \tbar_e) \fbar_e] 
    = \sum_{e \in E} \tbar_e f_e + \sum_{e \in E} \sup_{t_e \ge \tbar_e} (t_e - \tbar_e) (f_e - \fbar_e) 
    \longrightarrow \min_{f = \Theta x: x \in X} \\
    &\Longleftrightarrow \Psi(f) = \sum_{e \in E} f_e \tbar_e 
    \longrightarrow \min_{\substack{f = \Theta x : \\ x \in X,\, f_e \le \fbar_e}},\tag{SD}\label{PrimalSD}
\end{align*}
and its dual problem is
\begin{align*}\label{DualSD}
    Q(t) &= - \inf_{f = \Theta x : x \in X} \sum_{e \in E} [t_e f_e - (t_e - \tbar_e) \fbar_e]\\
    &= - \sum_{w \in OD} d_w T_w(t) + \langle t - \tbar, \fbar \rangle
    \longrightarrow \min_{t_e \ge \tbar_e}. \tag{DualSD}
\end{align*}
In contrast with the Beckmann model, the equilibrium state in the stable dynamics model is defined by pair $(f^*, t^*)$ (in particular, it differs from the system optimum $(f^*, \tbar)$ in the model only by the time value).

\section{Numerical methods}
\label{sec:methods}

We have the following objective functions
\begin{itemize}
    \item The stable dynamics model:
    \[
    Q(t) = \underbrace{- \sum_{w \in OD} d_w T_w(t)}_{\Phi(t)} + \underbrace{\langle t - \tbar, \fbar \rangle}_{h(t)},
    \]
    \item The Beckmann model:
    \[
    Q(t) = \underbrace{- \sum_{w \in OD} d_w T_w(t)}_{\Phi(t)} + \underbrace{\sum_{e \in E} \fbar_e \left(\frac{t_e - \tbar_e}{\tbar_e \rho} \right)^{\mu} \frac{\left (t_e - \bar{t}_e \right)}{1 + \mu}}_{h(t)}.
    \]
\end{itemize}
In both cases it has form
\begin{equation}
\label{eq::composite_Q}
    Q(t) = \Phi(t) + h(t) \longrightarrow \min_{t \ge \bar{t}}.  
\end{equation}
The optimization problem~\eqref{eq::composite_Q} is convex, non-smooth and composite. We use all these properties to identify the best optimization method to solve the considered problem.

\subsection{Subgradient}
\label{subsec:subgrad}

In our research, we consider first-order methods, i.e.\ they require a subgradient of $\Phi(t)$, the properties and effective computation of which we discuss in this section.

To get the subdifferential $\sgr \Phi(t)$ let us re-write $\Phi(t)$ in the following way:
\[
\Phi(t) = -\sum_{w \in OD} d_{w} T_{w}(t) 
= -\sum_{w \in OD} d_{w} \min_{p \in P_w} \langle t, a_p\rangle ,
\]
where the vector $a_p = (\delta_{e p})_{e \in E}$ encodes a path $p$.
Obviously, the shortest path may not be unique.
Using the rules of subgradient calculus \citep{rockafellar2015convex} we get the following expression: 
\[
    \sgr \Phi(t) = - \sum_{w \in OD} d_{w} \sgr \left(\min_{p \in P_w} \langle t, a_p \rangle\right) 
    = - \sum_{w \in OD} d_{w} \mathrm{Conv} \{a_p : p \in P_w,\; T_p(t) = T_w(t)\},
\]
i.e.\ the subdifferential $\sgr \Phi(t)$ is a sum of convex hulls of binary vectors that encode the shortest length paths. An important consequence is that for any $t_1, t_2 \in \R_+^{|E|}$ and $\nabla \Phi(t_1) \in \sgr \Phi(t_1)$, $\nabla \Phi(t_2) \in \sgr \Phi(t_2)$ the following bound holds:
\begin{equation}
\label{eq::M_const}
    \norm*{\nabla \Phi(t_1) - \nabla \Phi(t_2)}_2 \le M 
    = \sqrt{2 H} \sum_{w \in OD} d_{w},
\end{equation}
where $H$ is the diameter of the graph $G$.

Note that any element from the set $\sgr \Phi(t)$ has form $\nabla \Phi(t) = - f$, where $f = \Theta x$ is a flow distribution on links induced by $x \in X$ concentrated on the shortest paths for given times $t$ (and vice versa: any such $f$ corresponds to a subgradient of $\Phi(t)$).

In practice, the calculation of flows $f$ is the most expensive part, since we have to find the shortest paths for all pairs $w \in OD$. We use the following algorithm~\ref{alg::flows_aggregation}. We use Dijkstra's algorithm \citep{dijkstra1959note} to find the shortest paths in line~3, which runs in $O(|E| + |V| \log |V|)$ time; finding the traversal order with topological sort (Section~22.4 in \cite{cormen2009introduction}) and further flows aggregation have linear performance $O(|V|)$. Hence, the total complexity of Algorithm~\ref{alg::flows_aggregation} is $O\bigl(|O| (|E| + |V| \log |V|)\bigr)$. When the transportation network is an (almost) planar graph or another sparse graph, $|E| = O(|V|)$ and the complexity is $O(|O| \cdot |V| \log |V|)$. Moreover, flows reconstruction for every source $o \in O$ can be computed in parallel. And Dijkstra's algorithm also can be parallelized and has efficient implementations \citep{crauser1998parallelization, chao2010developed}.



\begin{algorithm}[!h]
    \caption{Flows reconstruction}
    \label{alg::flows_aggregation}
    \begin{algorithmic}[1]
    \REQUIRE times $t$
    \STATE $f \coloneqq \bm{0}_{|E|}$ \COMMENT{flows on edges}
    \FOR{origin $o$ in $O$}
        \STATE Get a shortest-path tree $\mathcal{T}_o$ from $o$ to all destinations in $D$ with weights $t$
        \STATE $\mathrm{traversal\_order} \coloneqq \mathrm{TopologicalSort}(\mathcal{T}_o)$
        \COMMENT{sorting from furthest to closest vertices}
        \STATE $f_{\mathrm{out}} \coloneqq \bm{0}_{|V|}$ \COMMENT{total output flow from each vertex}
        \STATE $f_{\mathrm{out}}[v] \coloneqq d_{w}$ for $w = (o, v) \in OD$
        \FOR{$v$ in $\mathrm{traversal\_order}$}
            \STATE Get predecessor $p$ of $v$ in $\mathcal{T}_o$
            \STATE $e \coloneqq (p, v)$
            \STATE $f[e] \coloneqq f[e] + f_{\mathrm{out}}[v]$ 
            \STATE $f_{\mathrm{out}}[p] \coloneqq f_{\mathrm{out}}[p] + f_{\mathrm{out}}[v]$
        \ENDFOR
    \ENDFOR
    \RETURN flows $f$
    \end{algorithmic}
\end{algorithm}

\subsection{Reconstruction of admissible flows in SD model}

For given times $t$ considered Algorithm~\ref{alg::flows_aggregation} reconstructs feasible flows $f$, i.e.\ $f = \Theta x$ for some $x \in X$.
These flows meet all the constraints in the Beckmann model, but they can violate the capacity constraints in the stable dynamics model. In the latter case, an additional step is required to obtain admissible flows from $f$. Note that we could instead find flows that meet capacity constraints first (Theorem~8 from \cite{nesterov2003stationary}), but to reconstruct feasible flows from them is a more complex problem.

Suppose we are given some flows $g = \Theta x$ such that 
\begin{equation}\label{eq:xi}
    \xi = 1 - \max_{e \in E} g_e / \fbar_e  > 0.
\end{equation}
Then for any $f = \Theta x$ we can construct admissible flows $\pi(f)$ in the following way: let $\eta = \max_{e \in E} f_e / \fbar_e - 1$, then
\[
\pi(f) = \begin{cases}
    f, & \eta \le 0, \\
    \frac{\xi f + \eta g}{\xi + \eta}, & \eta > 0.
\end{cases}
\]

In practice, we propose the following procedure to find admissible flows $g$: run some optimization method (e.g.\ UGM) for a small number of iterations for the same problem but with decreased capacities:
$\frac{1}{2} \fbar$ instead of $\fbar$; if obtained feasible flows $\hat{f}^N$ satisfy $\hat{f}^N \le \frac{3}{4} \fbar$, then take $g = \hat{f}^N$; otherwise, run it again with capacities $\frac{3}{4} \fbar$ and check $\hat{f}^N \le \frac{7}{8} \fbar$, etc.

\paragraph{Stopping criterion.}
The stopping criterion we use for the stable dynamics model is based on a duality gap
\begin{equation}\label{eq:UGM_stable_stop}
    Q(\hat{t}^N) + \Psi(\pi(\hat{f}^N)) \le \eps,
\end{equation}
where $\hat{f}^N \in \left\{\Theta x : x \in X\right\}$, $\hat{t}^N \ge \bar{t}$ are estimates of an equilibrium $(f^*, t^*)$ after $N$ iterations of the applied method. Note that here the duality gap with $\hat{f}^N$ is not applicable.


\subsection{Universal gradient method}
\label{subsec:ugd}

The method for solving non-smooth problems with smooth techniques was proposed by~\cite{nesterov2015universal} and was called \emph{the universal gradient method}.
The pseudocode of UGM for the considered problem~\eqref{eq::composite_Q} is provided in Algorithm~\ref{alg::univ_gd}. Here the euclidean prox-structure is used. Note that we did not specify the stopping criterion as it can be different for different models.

\begin{algorithm}[!h]
    \caption{Universal gradient method}
    \label{alg::univ_gd}
    \begin{algorithmic}[1]
    \REQUIRE $L_0 > 0$, accuracy $\eps > 0$
    \STATE Set $t^0 \coloneqq \tbar$, $k \coloneqq 0$
    \REPEAT
        \STATE $L_{k+1} \coloneqq L_k / 2$ 
        \WHILE{\TRUE}
            \STATE $t^{k+1} \coloneqq \argmin\limits_{t \in \dom h} \langle \nabla \Phi(t^k), t - t^k \rangle + h(t) + L_{k+1} \frac{\norm*{t - t^k}_2^2}{2}$
            \IF{$\Phi(t^{k+1}) \le \Phi(t^{k}) + \left\langle \nabla \Phi(t^{k}), t^{k+1} - t^{k} \right\rangle + L_{k+1} \frac{\norm*{t^{k+1} - t^k}_2^2}{2} + \frac{\eps}{2}$}
                \STATE \BREAK
            \ELSE
                \STATE $L_{k+1} \coloneqq 2 L_{k+1} $
            \ENDIF
        \ENDWHILE
        \STATE $k \coloneqq k + 1$
    \UNTIL{Stopping criterion is fulfilled}
    \end{algorithmic}
\end{algorithm}

Now let us define 
\begin{equation}
\label{eq:reconstruct_UGM}
    \hat{f}^N = - \frac{1}{S_N} \sum_{k = 0}^{N - 1} \frac{\nabla \Phi(t^k)}{L_{k+1}}, \quad
    \hat{t}^N = \frac{1}{S_N} \sum_{k = 1}^{N} \frac{t^k}{L_k}, \quad
    S_N = \sum_{k = 1}^{N} \frac{1}{L_k},
\end{equation}
where $L_k$ are the estimates of the local Lipschitz constant in UGM and UMST methods.

Convergence of the UGM was proved in~\cite{nesterov2015universal} and is summarized in the following lemma and theorem. 

\begin{lemma}\label{lem:UGM_stable}
    After $N$ iterations of UGM for the stable dynamics model it holds that
    \begin{gather}
        Q(\hat{t}^N) - Q(t^*) \le \frac{R^2}{S_N} + \frac{\eps}{2}, \label{eq:UGM_stable_Q}\\
        0 \le Q(\hat{t}^N) + \Psi(\hat{f}^N) + \langle t^* - \tbar, (\hat{f}^N - \fbar)_+ \rangle \le \frac{R^2}{S_N} + \frac{\eps}{2}, \label{eq:UGM_stable_dg}\\
        \norm{(\hat{f}^N - \fbar)_+}_2 \le \frac{4 R}{S_N} + \sqrt{\frac{2 \eps}{S_N}}, \label{eq:UGM_stable_f}
    \end{gather}
    where $\hat{f}^N$, $\hat{t}^N$, and $S_N$ are defined by~\eqref{eq:reconstruct_UGM}, and $R = \norm{t^* - \tbar}_2$ is the distance from the starting point to a solution.
\end{lemma}

\begin{theorem}\label{thm:UGM_stable}
    Let $L_0 \le \frac{M^2}{\eps}$, where $M$ comes from~\eqref{eq::M_const}.
    Then after at most
    \begin{equation}\label{eq:UGM_stable_complexity_Q}
        N_Q = 2 \left(\frac{R M}{\eps}\right)^2
    \end{equation}
    iterations of UGM for the stable dynamics model it holds that $Q(\hat{t}^N) - Q(t^*) \le \eps$.
    Moreover, the stopping criterion~\eqref{eq:UGM_stable_stop} is fulfilled after at most
    \begin{equation}\label{eq:UGM_stable_complexity_stop}
        N_{stop} 
        = O\left(\left(\frac{R M}{\eps}\right)^2 \max\left\{1,\,
        \left(\frac{\langle g - f^*, \tbar\rangle}{\xi R \min_e \fbar_e}\right)^2 \right\}\right)
    \end{equation}
    iterations, where $\xi$ comes from~\eqref{eq:xi}.
\end{theorem}

Now we provide results on the rate of convergence for the Beckmann model. The stopping criterion in this case is the following:
\begin{equation}\label{eq:UGM_Beckmann_stop}
    Q(\hat{t}^N) + \Psi(\hat{f}^N) \le \eps.
\end{equation}

\begin{lemma}\label{lem:UGM_Beckmann}
    After $N$ iterations of UGM for the Beckmann model it holds that
    \begin{gather*}\label{eq:UGM_Beckmann}
        Q(\hat{t}^N) - Q(t^*) \le \frac{R^2}{S_N} + \frac{\eps}{2},\\
        0 \le Q(\hat{t}^N) + \Psi(\hat{f}^N) \le \frac{\norm{\tau(\hat{f}^N) - \tbar}_2^2}{S_N} + \frac{\eps}{2},
    \end{gather*}
    where $\hat{f}^N$, $\hat{t}^N$, $S_N$ are defined by~\eqref{eq:reconstruct_UGM}, and $R = \norm{t^* - \tbar}_2$
\end{lemma}

\begin{theorem}\label{thm:UGM_Beckmann}
    Let $L_0 \le \frac{M^2}{\eps}$, where $M$ comes from~\eqref{eq::M_const}.
    Then after at most
    \begin{equation}\label{eq:UGM_Beckmann_complexity_Q}
        N_Q = 2 \left(\frac{R M}{\eps}\right)^2
    \end{equation}
    iterations of UGM for the Beckmann model it holds that $Q(\hat{t}^N) - Q(t^*) \le \eps$.
    Moreover, the stopping criterion~\eqref{eq:UGM_Beckmann_stop} is fulfilled after at most
    \begin{equation}\label{eq:UGM_Beckmann_complexity_stop}
        N_{stop} 
        = 2 \left(\frac{\tilde{R} M}{\eps}\right)^2
    \end{equation}
    iterations, where
    \begin{equation}\label{def:tilde_R}
        \tilde{R}^2 = \rho^2 \sum_{e \in E} \frac{\tbar_e^2}{\fbar_e^{2 / \mu}} \left(\sum_{w \in OD} d_w\right)^{2 / \mu}.
    \end{equation}
\end{theorem}

\subsection{Universal Method of Similar Triangles}
\label{subsec:ustm}

Let us introduce the following notations:
\[
\phi_0(t) = \frac{1}{2} \norm*{t - t^0}_2^2,
\]
\[
\phi_{k+1}(t) = \phi_k(t) + \alpha_{k+1} \left[\Phi(y^{k+1}) + \left\langle \nabla \Phi(y^{k+1}), t - y^{k+1} \right\rangle + h(t) \right].
\]

\begin{algorithm}[!h]
    \caption{Universal Method of Similar Triangles}
    \label{alg::univ_triangles}
    \begin{algorithmic}[1]
    \REQUIRE $L_0 > 0$, accuracy $\eps > 0$
    \STATE $u^0 = t^0 \coloneqq \tbar$, $A_0 \coloneqq 0$, $k \coloneqq 0$
    \REPEAT
        \STATE $L_{k+1} \coloneqq L_k / 2$ 
        \WHILE{\TRUE}
            \STATE 
            $\begin{cases} 
                \alpha_{k+1} \coloneqq \frac{1}{2 L_{k+1} } + \sqrt{\frac{1}{4 L_{k+1}^2} + \frac{A_k}{L_{k+1}} }, \quad A_{k+1} \coloneqq A_k + \alpha_{k+1}
                \\
                y^{k+1} \coloneqq \frac{\alpha_{k+1} u^k + A_k t^k}{A_{k+1}}, \quad u^{k+1} \coloneqq \argmin\limits_{t \in \dom h} \phi_{k+1}(t)
                \\
                t^{k+1} \coloneqq \frac{\alpha_{k+1} u^{k+1} + A_k t^k}{A_{k+1}}
            \end{cases}$
            \IF{$\Phi(t^{k+1}) \le \Phi(y^{k+1}) + \left\langle \nabla \Phi(y^{k+1}), t^{k+1} - y^{k+1} \right\rangle 
                + \frac{L_{k+1}}{2} \norm*{t^{k+1} - y^{k+1}}_2^2 + \frac{\alpha_{k+1}}{2 A_{k+1}} \eps$}
                \STATE \BREAK
            \ELSE
                \STATE $L_{k+1} \coloneqq 2 L_{k+1}$
            \ENDIF
        \ENDWHILE
        \STATE $k \coloneqq k + 1$
    \UNTIL{Stopping criterion is fulfilled}
    \end{algorithmic}
\end{algorithm}

Flows are reconstructed in the following way:
\begin{equation}
\label{eq:reconstruct_UMST}
    \hat{f}^N = - \frac{1}{A_N} \sum_{k = 1}^{N} \alpha_k \nabla \Phi(y^k)
\end{equation}

\begin{lemma}\label{lem:UMST_stable}
    After $N$ iterations of UMST for the stable dynamics model it holds that
    \begin{gather}
        Q(t^N) - Q(t^*) \le \frac{R^2}{A_N} + \frac{\eps}{2}, \label{eq:UMST_stable_Q} \\
        0 \le Q(t^N) + \Psi(\hat{f}^N) + \langle t^* - \tbar, (\hat{f}^N - \fbar)_+ \rangle \le \frac{R^2}{A_N} + \frac{\eps}{2}, \label{eq:UMST_stable_dg} \\ 
        \norm{(\hat{f}^N - \fbar)_+}_2 \le \frac{4 R}{A_N} + \sqrt{\frac{2 \eps}{A_N}}, \label{eq:UMST_stable_f}
    \end{gather}
    where $\hat{f}^N$ is defined by~\eqref{eq:reconstruct_UMST} and $R = \norm{t^* - \tbar}_2$ is the distance from the starting point to a solution.
\end{lemma}

\begin{theorem}\label{thm:UMST_stable}
    Let $L_0 \le \frac{4 M^2}{\eps}$, where $M$ comes from~\eqref{eq::M_const}.
    Then after at most
    \begin{equation}\label{eq:UMST_stable_complexity_Q}
        N_Q = 4 \left(\frac{R M}{\eps}\right)^2
    \end{equation}
    iterations of UGM for the stable dynamics model it holds that $Q(t^N) - Q(t^*) \le \eps$.
    Moreover, the stopping criterion~\eqref{eq:UGM_stable_stop} with $\hat{t}^N = t^N$ is fulfilled after at most
    \begin{equation}\label{eq:UMST_stable_complexity_stop}
        N_{stop}
        = O\left(\left(\frac{R M}{\eps}\right)^2 \max\left\{1,\,
        \left(\frac{\langle g - f^*, \tbar\rangle}{\xi R \min_e \fbar_e}\right)^2 \right\}\right)
    \end{equation}
    iterations, where $\xi$ comes from~\eqref{eq:xi}.
\end{theorem}

\begin{theorem}\label{thm:UMST_Beckmann}
    Let $L_0 \le \frac{4 M^2}{\eps}$, where $M$ comes from~\eqref{eq::M_const}.
    Then after at most
    \begin{equation}\label{eq:UMST_Beckmann_complexity_Q}
        N_Q = 4 \left(\frac{R M}{\eps}\right)^2
    \end{equation}
    iterations of UMST for the Beckmann model it holds that $Q(t^N) - Q(t^*) \le \eps$.
    Moreover, the stopping criterion~\eqref{eq:UGM_Beckmann_stop} with $\hat{t}^N = t^N$ is fulfilled after at most
    \begin{equation}\label{eq:UMST_Beckmann_complexity_stop}
        N_{stop} 
        = 4 \left(\frac{\tilde{R} M}{\eps}\right)^2
    \end{equation}
    iterations, where $\tilde{R}$ is defined by~\eqref{def:tilde_R}.
\end{theorem}

\subsection{Method of Weighted Dual Averages}
\label{subsec:pdsm}
%


\begin{algorithm}[!h]
    \caption{Method of Weighted Dual Averages}
    \label{alg::cpdsm}
    \begin{algorithmic}[1]
    \REQUIRE accuracy $\eps > 0$, constant $\chi > 0$
    \STATE $s^0 \coloneqq \vec{0}$, $t^0 \coloneqq \tbar$, $k \coloneqq 0$
    \REPEAT
        \STATE Compute subgradient $g^k$, set $s^{k + 1} \coloneqq s^{k} + \frac{1}{\norm{g^k}_2} g^k$
            \begin{itemize}
                \item non-composite case: $g^k \coloneqq \nabla \Phi(t^k) + \nabla h(t^k)$
                \item composite case:  $g^k \coloneqq \nabla \Phi(t^k)$
            \end{itemize}
        \STATE Set $\beta_{k+1} \coloneqq \frac{\hat{\beta}_{k+1}}{\chi}$, where $\hat{\beta}_{k+1} = \sum_{i = 0}^{k} \frac{1}{\hat{\beta}_{i}},\; \hat{\beta}_{0} = 1$
        \STATE Set $t^{k+1}$
            \begin{itemize}
                \item non-composite case: $t^{k+1} \coloneqq \argmin\limits_{t \in \dom h} \langle s^{k+1}, t \rangle + \frac{\beta_{k+1}}{2} \norm*{t - t^0}_2^2$
                \item composite case: $t^{k+1} \coloneqq \argmin\limits_{t \in \dom h} \langle s^{k+1}, t \rangle + \frac{\beta_{k+1}}{2} \norm*{t - t^0}_2^2 + \sum_{i = 0}^{k} \frac{1}{\norm{g^k}_2}h(t)$
            \end{itemize}
        \STATE $k \coloneqq k + 1$
    \UNTIL{Stopping criterion is fulfilled}
    \end{algorithmic}
\end{algorithm}

Convergence of WDA-method was proved in~\cite{nesterov2009primal} and is summarized in the following theorem. 

\begin{theorem}\label{thm:WDA}
    Non-composite WDA-method satisfies the following bounds
    \begin{itemize}
        \item for the stable dynamics model:
        \[
        Q(\hat{t}^k) - Q(t^*) = O\left(\frac{M + \norm{\fbar}_2}{\sqrt{k}} \left(\frac{R^2}{\chi} + \chi\right)\right),
        \]
        \item for the Beckmann model if $\mu \le 1$:
        \[
        Q(\hat{t}^k) - Q(t^*) = O\left(\frac{1}{\sqrt{k}} \left(M + \max_e \fbar_e \left[\frac{2 R + \chi}{\tbar_e \rho}\right]^\mu\right) \left(\frac{R^2}{\chi} + \chi\right)\right).
        \]
    \end{itemize}
\end{theorem}

\section{Numerical experiments}
\label{sec:experiments}

\begin{figure}[H]
    \centering
    \includegraphics[scale = 0.6, trim={0cm 0 0cm 0},clip]{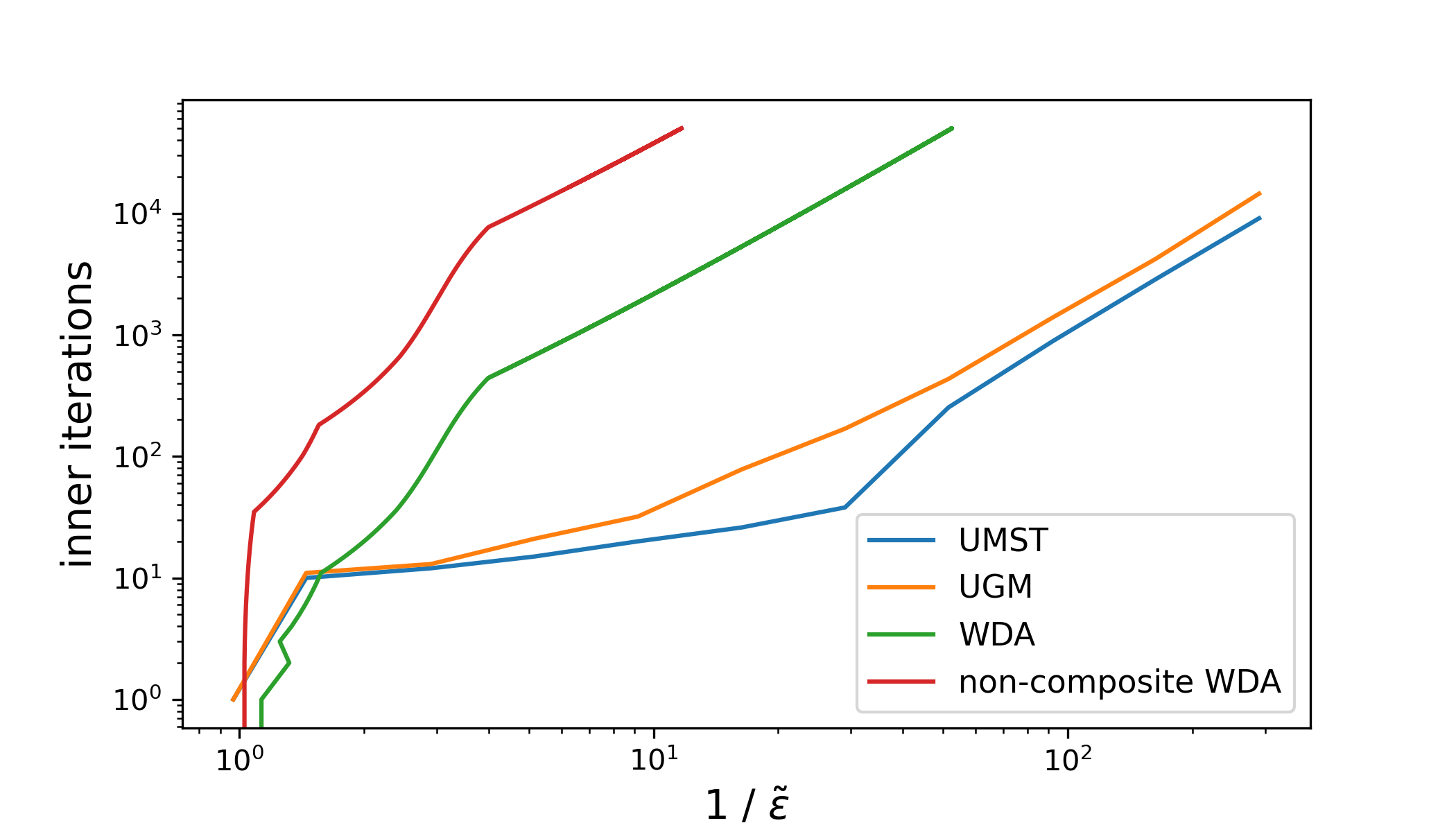}
    \caption{Convergence rates of UMST, UGM, composite and non-composite WDA-methods for the stable dynamics model with the stopping criterion~\eqref{eq:UGM_stable_stop}. Here $\tilde{\eps}$ is the relative accuracy $\eps / \Delta_0$, where $\Delta_0$ is the duality gap at the start point.}
    \label{fig:SD_convergence}
\end{figure}

This section presents numerical results for the algorithms described above, namely, composite variants of UMST and UGM, both composite and non-composite WDA-method, on the Anaheim network \citep{bstabler, chudak2007static}.
The network consists of 38 zones, 416 nodes, and 916 links. 
Experiments and the source code in Python 3\citep{python3} can be found in \cite{meruzakub}. We used Dijkstra's algorithm for finding the shortest paths in the network from the graph-tool library \citep{peixoto-graph-tool-2014}, where it is implemented in C++. We also used the Numpy library \citep{numpy} for all vector operations.
\paragraph{Stable dynamics model.} Parameters of the network are adjusted to the Beckmann model, so we have to increase the capacities to ensure the existence of an equilibrium for the stable dynamics model. In our experiments, the capacities are multiplied by $2.5$.
In Figure~\ref{fig:SD_convergence}, we plot the number of (inner) iterations of the algorithms required to fulfill the stopping criterion~\eqref{eq:UGM_stable_stop} against $1 / \eps$.
We consider the number of inner iterations for Alg.~\ref{alg::univ_triangles} and Alg.~\ref{alg::univ_gd} since the complexity of an inner iteration in this case is similar to the complexity of an iteration of the other algorithms. Note that according to \cite[formula (2.23)]{nesterov2015universal} the number $N(k)$ of inner iterations of UGM or UMST at step $k$ is bounded as
\[
N(k) \le 2 k + \log_2\left(\frac{M^2}{\eps L_0}\right),
\]
so asymptotic rates from Theorems~\ref{thm:UGM_stable}, \ref{thm:UGM_Beckmann}, \ref{thm:UMST_stable}, and \ref{thm:UMST_Beckmann} are still valid.

As we can see, the best results are shown by UMST, followed by UGM having similar performance.
Both composite and non-composite WDA-method are much slower.

\paragraph{Beckmann model.} For the Beckmann model we also compare our methods with the Frank--Wolfe algorithm (Alg.~\ref{alg::frank-wolfe})~--- which theoretical convergence rate for convex objective (with Lipschitz-continuous gradient) is $O(1/\eps)$ \citep{pedregosa2020linearly, jaggi2013revisiting}.
\begin{algorithm}[!h]
    \caption{Frank-Wolfe algorithm}
    \label{alg::frank-wolfe}
    \begin{algorithmic}[1]
    \REQUIRE accuracy $\eps > 0$
    \STATE $t^0 \coloneqq \tbar$, $f^0 \coloneqq \argmin\limits_{s \in \{\Theta x : x \in X\}} \langle t^0, s \rangle$, $k \coloneqq 0$
    \REPEAT
        \STATE $s^k \coloneqq \argmin\limits_{s \in \{\Theta x : x \in X\}} \langle t^k, s \rangle$, $t_e^k \coloneqq \frac{\partial \Psi (f^k)}{\partial f_e} = \tau_e(f^k)$ 
        \STATE $\gamma_k \coloneqq \frac{2}{k + 2}$, $f^{k + 1} \coloneqq (1 - \gamma_k) f^k + \gamma_k s^k$ 
        \STATE $k \coloneqq k + 1$
    \UNTIL{Stopping criterion is fulfilled}
    \end{algorithmic}
\end{algorithm}

\begin{figure}[H]
    \centering
    \includegraphics[scale = 0.6, trim={0cm 0 0cm 0},clip]{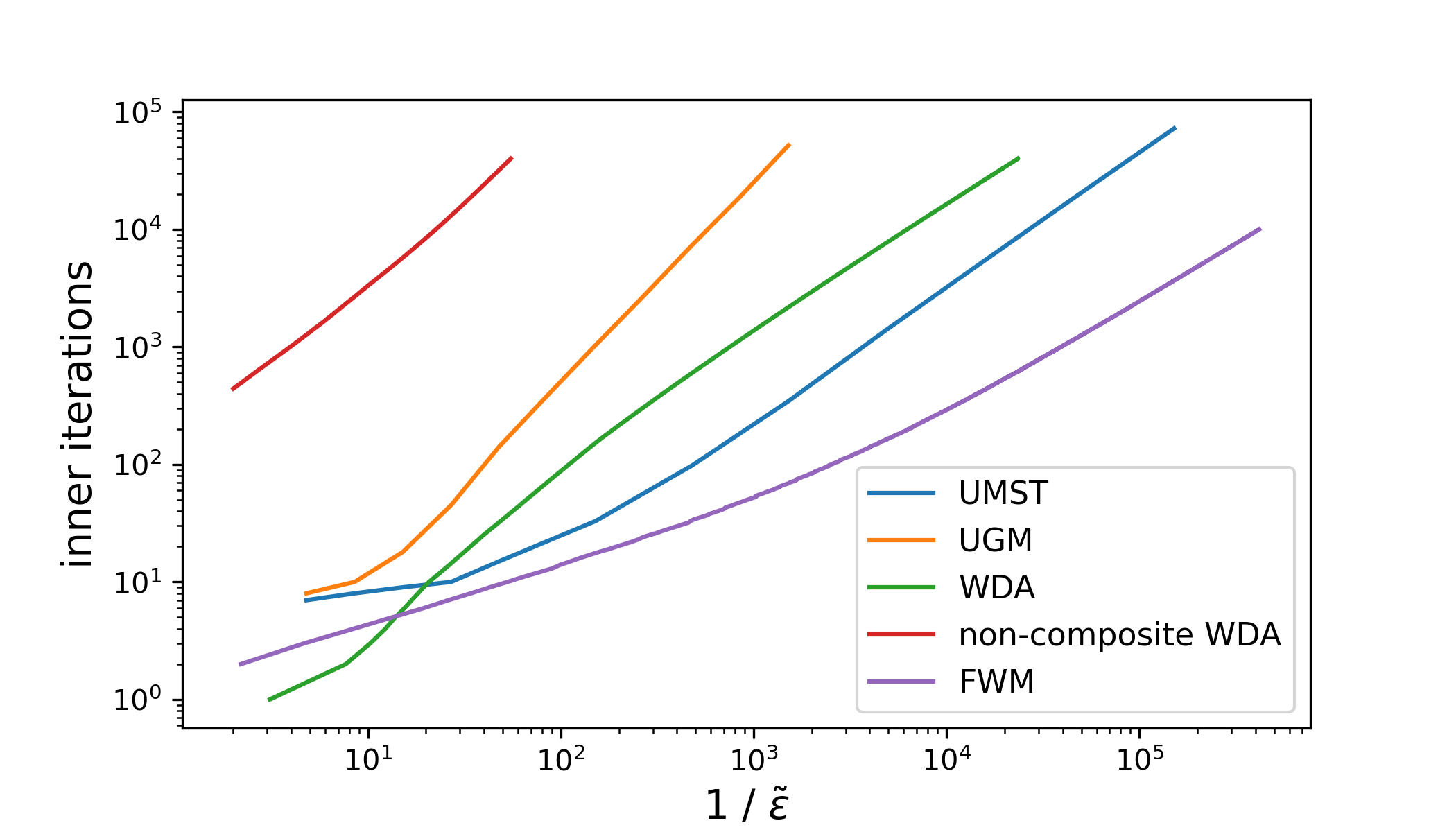}
    \caption{Convergence rates of UMST, UGM, composite and non-composite WDA-methods, and the Frank--Wolfe method for the Beckmann model with the stopping criterion~\eqref{eq:UGM_Beckmann_stop}. Here $\tilde{\eps}$ is the relative accuracy $\eps / \Delta_0$, where $\Delta_0$ is the duality gap at the start point.}
    \label{fig:Beckmann_convergence}
\end{figure}

Figure~\ref{fig:Beckmann_convergence} shows the convergence rates of the methods for the Beckmann model. The Frank--Wolfe method demonstrates the best results and is followed by UMST. Unlike the stable dynamics case, composite WDA-method is faster than UGM. However, the non-composite WDA-method has the worst performance again.
    
\section{Conclusion}
\label{sec:conclusion}

We considered several primal-dual subgradient methods for finding equilibria in the stable dynamics and the Beckmann models. We suggested a way to reconstruct admissible flows in the stable dynamics model, which provides us with a novel computable duality gap.
Complexity bounds for UMST and UGM were presented in terms of the iterations number required to achieve a desired accuracy in the dual function value or the duality gap.
Finally, we conducted numerical experiments comparing convergence of the considered algorithms on the Anaheim transportation network: UMST is the best one for optimization of the dual problems in both models. Also, using the duality gap as a stopping criterion, we compared these methods with the Frank--Wolfe algorithm for the Beckmann model~--- which, as expected, remains the most suitable approach in this case (but it is not applicable for the stable dynamics model).

The reader may be interested in another related topic, searching stochastic traffic equilibria. In \cite{meruza2018, baimurzina2019universal} we with our colleagues studied the application of the UMST for finding Nash--Wardrop stochastic equilibria in the Beckmann model. In this case, a driver selects a route randomly according to Gibbs’ distribution taking into account current time costs on the links of the network. 
It leads to iteration complexity $O(\frac{1}{\sqrt{\gamma \eps}})$, where $\gamma > 0$ is a stochasticity parameter (when $\gamma \to 0$ the model boils down to the ordinary Beckmann model). However, the great decrease in the number of iterations comes along with a more expensive calculation of the objective function's gradient.

\section*{Acknowledgements}
We would like to thanks Yu. Nesterov for fruitful discussions.

The research of M. Kubentayeva was supported by Russian Science Foundation (project 18-71-10108).

The research of A. Gasnikov was partially supported by RFBR, project number 18-29-03071 mk, and was partially supported by the Ministry of Science and Higher Education of the Russian Federation (Goszadaniye) no 075-00337-20-03.








\bibliography{lib}

@article{dijkstra1959note,
  title={A note on two problems in connexion with graphs},
  author={Dijkstra, Edsger W and others},
  journal={Numerische mathematik},
  volume={1},
  number={1},
  pages={269--271},
  year={1959}
}

@book{cormen2009introduction,
  title={Introduction to algorithms, 3rd-edition},
  author={Cormen, Thomas H and Leiserson, Charles E and Rivest, Ronald L and Stein, Clifford},
  year={2009},
  publisher={MIT press}
}

@inproceedings{crauser1998parallelization,
  title={A parallelization of {D}ijkstra's shortest path algorithm},
  author={Crauser, Andreas and Mehlhorn, Kurt and Meyer, Ulrich and Sanders, Peter},
  booktitle={International Symposium on Mathematical Foundations of Computer Science},
  pages={722--731},
  year={1998},
  organization={Springer}
}

@inproceedings{chao2010developed,
  title={Developed {D}ijkstra shortest path search algorithm and simulation},
  author={Chao, Yin and Hongxia, Wang},
  booktitle={2010 International Conference on Computer Design and Applications},
  volume={1},
  pages={V1--116},
  year={2010},
  organization={IEEE}
}

@book{bpr_functions,
  title={Traffic Assignment Manual},
  author={{US Bureau of Public Roads}},
  year={1964},
  publisher={Department of Commerce, Urban Planning Division, Washington D.C.}
}

@inproceedings{chudak2007static,
  title={Static traffic assignment problem: A comparison between {B}eckmann (1956) and {N}esterov \& de {P}alma (1998) models},
  author={Chudak, Fabian A and Dos Santos Eleuterio, Vania and Nesterov, Yurii},
  booktitle={7th Swiss Transport Research Conference},
  year={2007},
  organization={ETH}
}

@misc{bstabler,
  author = {{Transportation Networks for Research Core Team}},
  title = {Transportation Networks for Research},
  howpublished = {\url{https://github.com/bstabler/TransportationNetworks}},
  note = {Accessed: 2021-04-30},
  year = 2021
}

@misc{meruzakub,
  author = {Kubentayeva, Meruza},
  title = {Transport{N}et},
  howpublished = {\url{https://github.com/MeruzaKub/TransportNet}},
  note = {Accessed: 2021-04-30},
  year = 2021
}

@article{peixoto-graph-tool-2014,
         title = {The graph-tool python library},
         url = {\url{http://figshare.com/articles/graph_tool/1164194}},
         doi = {10.6084/m9.figshare.1164194},
         urldate = {2014-09-10},
         journal = {figshare},
         author = {Peixoto, Tiago P.},
         year = {2014}
}

@book{python3,
 author = {Van Rossum, Guido and Drake, Fred L.},
 title = {Python 3 Reference Manual},
 year = {2009},
 isbn = {1441412697},
 publisher = {CreateSpace},
 address = {Scotts Valley, CA}
}

@article{numpy,
 title         = {Array programming with {NumPy}},
 author        = {Charles R. Harris and K. Jarrod Millman and St{\'{e}}fan J. van der Walt and Ralf Gommers and Pauli Virtanen and David Cournapeau and Eric Wieser and Julian Taylor and Sebastian Berg and Nathaniel J. Smith and Robert Kern and Matti Picus and Stephan Hoyer and Marten H. van Kerkwijk and Matthew Brett and Allan Haldane and Jaime Fern{\'{a}}ndez del R{\'{i}}o and Mark Wiebe and Pearu Peterson and Pierre G{\'{e}}rard-Marchant and Kevin Sheppard and Tyler Reddy and Warren Weckesser and Hameer Abbasi and Christoph Gohlke and Travis E. Oliphant},
 year          = {2020},
 month         = sep,
 journal       = {Nature},
 volume        = {585},
 number        = {7825},
 pages         = {357--362},
 doi           = {10.1038/s41586-020-2649-2},
 publisher     = {Springer Science and Business Media {LLC}},
 url           = {https://doi.org/10.1038/s41586-020-2649-2}
}

@article{baimurzina2019universal,
  title={Universal method of searching for equilibria and stochastic equilibria in transportation networks},
  author={Baimurzina, D R and Gasnikov, A V and Gasnikova, E V and Dvurechensky, P E and Ershov, E I and Kubentaeva, M B and Lagunovskaya, A A},
  journal={Computational Mathematics and Mathematical Physics},
  volume={59},
  number={1},
  pages={19--33},
  year={2019},
  publisher={Springer}
}

@article{meruza2018,
  title={Searching stochastic equilibria in transport networks by universal primal-dual gradient method},
  author={Gasnikov, Alexander Vladimirovich and Kubentayeva, Meruza B},
  journal={Computer research and modeling},
  volume={10},
  number={3},
  pages={335--345},
  year={2018},
  publisher={Institute of Computer Science}
}

@techreport{beckmann1956studies,
  title={Studies in the Economics of Transportation},
  author={Beckmann, Martin J and McGuire, Charles B and Winsten, Christopher B},
  year={1956},
  publisher={Rand Corporation}
}

@book{patriksson2015traffic,
  title={The traffic assignment problem: models and methods},
  author={Patriksson, Michael},
  year={2015},
  publisher={Courier Dover Publications}
}

@book{rockafellar2015convex,
  title={Convex analysis},
  author={Rockafellar, Ralph Tyrell},
  year={2015},
  publisher={Princeton university press}
}

@book{sheffi1985urban,
  title={Urban transportation networks},
  author={Sheffi, Yosef},
  volume={6},
  year={1985},
  publisher={Prentice-Hall, Englewood Cliffs, NJ}
}

@article{nesterov2015universal,
  title={Universal gradient methods for convex optimization problems},
  author={Nesterov, Yu},
  journal={Mathematical Programming},
  volume={152},
  number={1-2},
  pages={381--404},
  year={2015},
  publisher={Springer}
}

@article{gasnikov2018universal,
  title={Universal method for stochastic composite optimization problems},
  author={Gasnikov, Alexander Vladimirovich and Nesterov, Yu E},
  journal={Computational Mathematics and Mathematical Physics},
  volume={58},
  number={1},
  pages={48--64},
  year={2018},
  publisher={Springer}
}

@article{nesterov2003stationary,
  title={Stationary dynamic solutions in congested transportation networks: summary and perspectives},
  author={Nesterov, Yurii and De Palma, Andre},
  journal={Networks and Spatial Economics},
  volume={3},
  number={3},
  pages={371--395},
  year={2003},
  publisher={Springer}
}

@article{nesterov2009primal,
  title={Primal-dual subgradient methods for convex problems},
  author={Nesterov, Yurii},
  journal={Mathematical programming},
  volume={120},
  number={1},
  pages={221--259},
  year={2009},
  publisher={Springer}
}

@article{frank1956algorithm,
  title={An algorithm for quadratic programming},
  author={Frank, Marguerite and Wolfe, Philip},
  journal={Naval research logistics quarterly},
  volume={3},
  number={1-2},
  pages={95--110},
  year={1956},
  publisher={Wiley Online Library}
}

@inproceedings{pedregosa2020linearly,
  title={Linearly Convergent Frank-Wolfe with Backtracking Line-Search},
  author={Pedregosa, Fabian and Negiar, Geoffrey and Askari, Armin and Jaggi, Martin},
  booktitle={Proceedings of the 23rdInternational Conference on Artificial Intelligence and Statistics},
  year={2020},
  url={https://arxiv.org/pdf/1806.05123.pdf}
}

@inproceedings{jaggi2013revisiting,
  title={Revisiting {F}rank--{W}olfe: Projection-free sparse convex optimization.},
  author={Jaggi, Martin},
  booktitle={Proceedings of the 30th international conference on machine learning},
  pages={427--435},
  year={2013}
}

@article{leblanc1985improved,
  title={Improved efficiency of the {F}rank--{W}olfe algorithm for convex network programs},
  author={LeBlanc, Larry J and Helgason, Richard V and Boyce, David E},
  journal={Transportation Science},
  volume={19},
  number={4},
  pages={445--462},
  year={1985},
  publisher={INFORMS}
}

@article{arezki1990full,
  title={A full analytical implementation of the {PARTAN}/{F}rank--{W}olfe algorithm for equilibrium assignment},
  author={Arezki, Y and Van Vliet, D},
  journal={Transportation Science},
  volume={24},
  number={1},
  pages={58--62},
  year={1990},
  publisher={INFORMS}
}

@article{fukushima1984modified,
  title={A modified {F}rank--{W}olfe algorithm for solving the traffic assignment problem},
  author={Fukushima, Masao},
  journal={Transportation Research Part B: Methodological},
  volume={18},
  number={2},
  pages={169--177},
  year={1984},
  publisher={Elsevier}
}

@article{chen2002faster,
author = {Chen, Anthony and Jayakrishnan, R. and Tsai, Wei},
year = {2002},
month = {01},
pages = {},
title = {Faster {F}rank--{W}olfe Traffic Assignment with New Flow Update Scheme},
volume = {128},
journal = {Journal of Transportation Engineering-asce - J TRANSP ENG-ASCE},
doi = {10.1061/(ASCE)0733-947X(2002)128:1(31)}
}

\section{Appendix}
\label{sec:appendix}

\subsection{Proofs for UGM}
\begin{proof}[\textbf{Proof of Lemma \ref{lem:UGM_stable}}]
    Note that function $\Phi(t)$ satisfies~\eqref{eq::M_const}.
    Then according to Theorem~1 in \cite{nesterov2015universal} applied with $\nu = 0$ one has
    \begin{align}\label{eq:Nesterov2015_inequality}
        Q(\hat{t}^N) 
        &\le \frac{1}{S_N} \sum_{k = 1}^{N} \frac{1}{L_k} Q(t^k) \nonumber\\
        &\le \min_{t \ge \tbar} \left\{\frac{1}{S_N} \sum_{k = 0}^{N - 1} \frac{1}{L_{k+1}} \left[\Phi(t^k) + \langle \nabla \Phi(t^k), t - t^k\rangle \right] + h(t) + \frac{\norm{t - t^0}_2^2}{S_N}\right\} + \frac{\eps}{2}.
    \end{align}
    Equation~\eqref{eq:UGM_stable_Q} follows immediately if one substitutes $t = t^*$.
    Now let us estimate the first term on the r.h.s.
    \begin{align*}
        \min_{t \ge \tbar} &\left\{\frac{1}{S_N} \sum_{k = 0}^{N - 1} \frac{1}{L_{k+1}} \left[\Phi(t^k) + \langle \nabla \Phi(t^k), t - t^k \rangle\right] + h(t) + \frac{\norm{t - \tbar}_2^2}{S_N}\right\} \\
        &= \min_{t \ge \tbar} \left\{\frac{1}{S_N} \sum_{k = 0}^{N - 1} \frac{1}{L_{k+1}} \underbrace{\left[\Phi(t^k) + \langle \nabla \Phi(t^k), 0 - t^k\rangle \right]}_{\le \Phi(0)} - \langle \hat{f}^N, t\rangle + \langle \fbar, t - \tbar\rangle + \frac{\norm{t - \tbar}_2^2}{S_N} \right\} \\
        &\le \Phi(0) - \langle \hat{f}^N, \tbar\rangle + \min_{t \ge \tbar} \left\{\langle \fbar - \hat{f}^N, t - \tbar\rangle + \frac{\norm{t - \tbar}_2^2}{S_N}\right\} \\
        &= - \Psi(\hat{f}^N) - \frac{S_N \norm{(\hat{f}^N - \fbar)_+}_2^2}{4}.
    \end{align*}
    Here we used that $\Phi(0) = - \sum_{w \in OD} d_w T_w(0) = 0$.
    Therefore,
    \[
        Q(\hat{t}^N) + \Psi(\hat{f}^N) + \frac{S_N \norm{(\hat{f}^N - \fbar)_+}_2^2}{4} \le \frac{\eps}{2}.
    \]
    Now notice that, since the flow $\hat{f}^N$ is induced by some traffic distribution $x \in X$, we have
    \begin{align*}
        0 &\le \Phi(t^*) + \langle t^*, \hat{f}^N\rangle \\
        & = Q(t^*) - \langle t^* - \tbar, \fbar \rangle + \Psi(\hat{f}^N) - \langle \tbar, \hat{f}^N\rangle + \langle t^*, \hat{f}^N\rangle\\
        & = Q(t^*) + \Psi(\hat{f}^N) + \langle t^* - \tbar, \hat{f}^N - \fbar \rangle \\
        & \le Q(\hat{t}^N) + \Psi(\hat{f}^N) + \langle t^* - \tbar, (\hat{f}^N - \fbar)_+ \rangle,
    \end{align*}
    hence
    \[
        Q(\hat{t}^N) + \Psi(\hat{f}^N) 
        \ge - \langle t^* - \tbar, (\hat{f}^N - \fbar)_+ \rangle
        \ge - R \norm{(\hat{f}^N - \fbar)_+}_2.
    \]
    This yields
    \[
    \frac{S_N \norm{(\hat{f}^N - \fbar)_+}_2^2}{4} - R \norm{(\hat{f}^N - \fbar)_+}_2 \le \frac{\eps}{2},
    \]
    and thus
    \[
        \norm{(\hat{f}^N - \fbar)_+}_2
        \le \frac{2 R}{S_N} \left(1 + \sqrt{1 + \frac{\eps S_N}{2 R^2}}\right) 
        \le \frac{4 R}{S_N} + \sqrt{\frac{2 \eps}{S_N}}.
    \]
\end{proof}

\begin{proof}[\textbf{Proof of Theorem \ref{thm:UGM_stable}}]
    Theorem~1 in \cite{nesterov2015universal} ensures that $L_k \le \frac{M^2}{\eps}$ for all $k \ge 0$, thus $S_N \ge \frac{\eps N}{M^2}$.
    Then the first bound~\eqref{eq:UGM_stable_complexity_Q} follows immediately from~\eqref{eq:UGM_stable_Q}.
    
    Now let us prove the second bound. First, suppose $\hat{f}^N_e \le \fbar_e$ for all $e \in E$. Then $\pi(\hat{f}^N) = \hat{f}^N$, thus by~\eqref{eq:UGM_stable_dg} for $N = N_Q$
    \[
        Q(\hat{t}^N) + \Psi(\pi(\hat{f}^N)) 
        = Q(\hat{t}^N) + \Psi(\hat{f}^N) 
        \le \frac{R^2}{S_N} + \frac{\eps}{2}
        \le \frac{(R M)^2}{\eps N} + \frac{\eps}{2}
        \le \eps.
    \]
    Otherwise, if $\hat{f}^N_e \not\le \fbar_e$, one has $\pi(\hat{f}^N) = \frac{\xi \hat{f}^N + \eta g}{\xi + \eta}$, where $\eta = \max_{e \in E} \hat{f}^N_e / \fbar_e - 1$, hence~\eqref{eq:UGM_stable_Q} and~\eqref{eq:UGM_stable_dg} yield
    \begin{align*}
        Q(\hat{t}^N) + \Psi(\pi(\hat{f}^N)) 
        & \le \frac{\xi}{\xi + \eta} \left(Q(\hat{t}^N) + \Psi(\hat{f}^N)\right) + \frac{\eta}{\xi + \eta} \left(Q(\hat{t}^N) + \Psi(g)\right) \\
        & = \frac{\xi}{\xi + \eta} \left(Q(\hat{t}^N) + \Psi(\hat{f}^N)\right) + \frac{\eta}{\xi + \eta} \left(Q(\hat{t}^N) - Q(t^*)\right) + \frac{\eta}{\xi + \eta} \left(\Psi(g) - \Psi(f^*)\right) \\
        & \le \frac{R^2}{S_N} + \frac{\eps}{2} + \frac{\eta}{\xi} \langle g - f^*, \tbar\rangle.
    \end{align*}
    Finally, according to~\eqref{eq:UGM_stable_f}
    \[
        \eta 
        = \max_{e \in E} \hat{f}^N_e / \fbar_e - 1 
        = \norm*{\frac{(\hat{f}^N - \fbar)_+}{\fbar}}_\infty
        \le \frac{1}{\min_e \fbar_e} \norm*{(\hat{f}^N - \fbar)_+}_2
        \le \frac{1}{\min_e \fbar_e} \left(\frac{4 R}{S_N} + \sqrt{\frac{2 \eps}{S_N}}\right).
    \]
    Combining all bounds together we obtain
    \[
    Q(\hat{t}^N) + \Psi(\pi(\hat{f}^N)) 
    \le \frac{R^2 M^2}{\eps N} + \frac{\langle g - f^*, \tbar\rangle}{\xi \min_e \fbar_e} \left(\frac{4 R M^2}{\eps N} + \sqrt{\frac{2 M^2}{N}}\right) + \frac{\eps}{2},
    \]
    and substituting $N = N_{stop}$, we conclude that the stopping criterion~\eqref{eq:UGM_stable_stop} is fulfilled.
\end{proof}

\begin{proof}[\textbf{Proof of Lemma \ref{lem:UGM_Beckmann}}]
    First of all, note that
    \[
        \max_{t \ge \tbar} \left\{\langle \hat{f}^N, t\rangle - \sum_{e \in E} \sigma^*_e (t_e) \right\} 
        = \sum_{e \in E} \sigma_e (\hat{f}^N_e) = \Psi(\hat{f}^N),
    \]
    and maximum is attained at point $t = \nabla \Psi(\hat{f}^N) = \tau(\hat{f}^N)$.
    As in the proof of Theorem~\ref{thm:UGM_stable}, the inequality~\eqref{eq:Nesterov2015_inequality} holds in Beckmann's model case. 
    Then the first term in the r.h.s.\ can be estimated as follows
    \begin{align*}
        \min_{t \ge \tbar} &\left\{\frac{1}{S_N} \sum_{k = 0}^{N - 1} \frac{1}{L_{k+1}} \left[\Phi(t^k) + \langle \nabla \Phi(t^k), t - t^k \rangle\right] + h(t) + \frac{\norm{t - \tbar}_2^2}{S_N}\right\} \\
        &= \min_{t \ge \tbar} \left\{\frac{1}{S_N} \sum_{k = 0}^{N - 1} \frac{1}{L_{k+1}} \underbrace{\left[\Phi(t^k) + \langle \nabla \Phi(t^k), 0 - t^k\rangle \right]}_{\le \Phi(0)} - \langle \hat{f}^N, t\rangle + 
        \sum_{e \in E} \sigma^*_e (t_e) + \frac{\norm{t - \tbar}_2^2}{S_N} \right\} \\
        &\le \Phi(0) + \left\{\sum_{e \in E} \sigma^*_e (t_e(\hat{f}^N_e)) - \langle \hat{f}^N, \tau(\hat{f}^N)\rangle + \frac{\norm{\tau(\hat{f}^N) - \tbar}_2^2}{S_N}\right\} \\
        & = - \Psi(\hat{f}^N) + \frac{1}{S_N} \norm{\tau(\hat{f}^N) - \tbar}_2^2,
    \end{align*}
    and we finally get an upper bound on the duality gap:
    \[
    0 \le Q(\hat{t}^N) + \Psi(\hat{f}^N) 
    \le \frac{\norm{\tau(\hat{f}^N) - \tbar}_2^2}{S_N} + \frac{\eps}{2}.
    \]
    In the same time, substituting $t = t^*$ one obtains
    \[
    Q(\hat{t}^N) \le Q(t^*) + \frac{\norm{t^* - \tbar}^2}{S_N} + \frac{\eps}{2}.
    \]
\end{proof}

\begin{proof}[\textbf{Proof of Theorem \ref{thm:UGM_Beckmann}}]
    By construction, $\hat{t}^N_e \le \sum_{w \in OD} d_w$ for all $e \in E$, thus $\norm{\tau(\hat{f}^N) - \tbar}_2 \le \tilde{R}$. According to Theorem~1 in \cite{nesterov2015universal} $S_N \ge \frac{\eps N}{M^2}$, thus the statement follows immediately from Lemma~\ref{lem:UGM_Beckmann}.
\end{proof}

\subsection{Proofs for UMST}
\begin{proof}[\textbf{Proof of Lemma \ref{lem:UMST_stable}}]
    According to the inequality~(30) in \cite{gasnikov2018universal}
    \begin{equation}
        Q(t^N) \le  \min_{t \ge \tbar} \left\{\frac{1}{A_N} \sum_{k = 1}^N \alpha_k \left[\Phi(y^k) + \left\langle \nabla \Phi(y^k), t - y^k \right\rangle\right] + h(t) + \frac{\norm{t - t^0}_2^2}{2 A_N}\right\} + \frac{\eps}{2}.
    \end{equation}
    Note that the above inequality has the same form as~\eqref{eq:Nesterov2015_inequality}, if one replaces $S_N$ with $A_N$, $\frac{1}{L_{k+1}}$ with $\alpha_k$, $y^k$ with $t^k$, and $\frac{\norm{t - t^0}_2^2}{S_N}$ with $\frac{\norm{t - t^0}_2^2}{2 A_N}$.
    Then the claim follows by the same reasoning as in the proof of Lemma~\ref{lem:UGM_stable}.
\end{proof}

\begin{proof}[\textbf{Proof of Theorem \ref{thm:UMST_stable}}]
    Due to \eqref{eq::M_const} one has
    \[
    \Phi(t^{k+1}) \le \Phi(y^{k+1}) + \langle \nabla \Phi(y^{k+1}), t^{k+1} - y^{k+1} \rangle + M \norm{t^{k+1} - y^{k+1}}_2.
    \]
    From Young's inequality we get that
    \[
    M \norm{t^{k+1} - y^{k+1}}_2 \le \frac{\alpha_{k+1}}{2 A_{k+1}} \eps + \frac{A_{k+1} M^2}{2 \alpha_{k+1} \eps} \norm{t^{k+1} - y^{k+1}}_2^2. 
    \]
    If $L_{k+1} \ge \frac{A_{k+1} M^2}{\alpha_{k+1} \eps}$, then the stopping condition for inner iterations is fulfilled. Therefore, at the end of the $k$-th iteration either $L_{k+1} < \frac{2 A_{k+1} M^2}{\alpha_{k+1} \eps}$ or $L_{k+1} = \frac{L_k}{2}$. 
    
    Now we are going to prove by induction that $\alpha_k \ge \frac{\eps}{2 M^2}$, which is equivalent to 
    $L_k \le \frac{2 M^2}{\eps} + \frac{4 M^4}{\eps^2} A_k$, for all $k \ge 1$. For $k = 1$ it follows from $A_1 = \alpha_1$ and $L_0 \le \frac{4 M^2}{\eps}$.
    In case where $L_{k+1} < \frac{2 A_{k+1} M^2}{\alpha_{k+1} \eps}$ equation $A_{k+1} = L_{k+1} \alpha_{k+1}^2$ immediately yields $\alpha_{k+1} \ge \frac{\eps}{2 M^2}$. If $L_{k+1} = \frac{L_k}{2}$, then by the induction hypothesis and monotonicity of the sequence $\{A_k\}_{k \in \N}$ we obtain
    \[
    L_{k+1} \le \frac{M^2}{\eps} + \frac{2 M^4}{\eps^2} A_{k-1} < \frac{2 M^2}{\eps} + \frac{4 M^4}{\eps^2} A_k.
    \]
    Therefore,
    \begin{equation}\label{eq:A_bound}
        A_N \ge \frac{\eps N}{2 M^2}.
    \end{equation}
    Arguing in the same way as in the proof of Theorem~\ref{thm:UGM_stable} we obtain that
    \[
    Q(t^N) - Q(t^*) \le \frac{2 R^2 M^2}{\eps N} + \frac{\eps}{2}
    \]
    and
    \[
    Q(\hat{t}^N) + \Psi(\pi(\hat{f}^N)) 
    \le \frac{2 R^2 M^2}{\eps N} + \frac{\langle g - f^*, \tbar\rangle}{\xi \min_e \fbar_e} \left(\frac{8 R M^2}{\eps N} + \sqrt{\frac{4 M^2}{N}}\right) + \frac{\eps}{2}.
    \]
    After substitution $N = N_Q$ or $N = N_{stop}$ the claim follows.
\end{proof}

\begin{proof}[\textbf{Proof of Theorem \ref{thm:UMST_Beckmann}}]
    Repeating the proof of Theorem~\ref{thm:UGM_Beckmann} we obtain that
    \[
    Q(t^N) - Q(t^*) \le \frac{R^2}{A_N} + \frac{\eps}{2}, \quad
    Q(t^N) + \Psi(\hat{f}^N) \le \frac{\tilde{R}^2}{A_N} + \frac{\eps}{2}.
    \]
    Then we conclude applying~\eqref{eq:A_bound}.
\end{proof}

\subsection{Proof of Theorem 3.8 (WDA)}
    According to eq.~(3.5) from \cite{nesterov2009primal},
    \[
    Q(\hat{t}^k) - Q(t^*) = O\left(\frac{L}{\sqrt{k}} \left(\frac{R^2}{\chi} + \chi\right)\right),
    \]
    whenever $\norm{g_k}_2 \le L$ for all $k$.
    
    In case of the stable dynamics model $\nabla h(t) = \fbar$, thus we can take $L = M + \norm{\fbar}_2$.
    
    For the Beckmann model 
    \[
    \frac{\partial h(t)}{\partial t_e} = \fbar_e \left(\frac{t_e - \tbar_e}{\tbar_e \rho} \right)^{\mu}.
    \]
    Theorem~3 in \cite{nesterov2009primal} yields that $\norm{t^k - t^*}_2^2 \le R^2 + \chi^2$ for all $k$, thus $\norm{t^k - \tbar}_2 \le 2 R + \chi$. Then using $\mu \le 1$ one obtains
    \[
    \norm{\nabla h(t^k)}_2 \le \norm{t^k - \tbar}_2^\mu \max_e \frac{\fbar_e}{(\tbar_e \rho)^\mu} 
    \le (2 R + \chi)^\mu \max_e \frac{\fbar_e}{(\tbar_e \rho)^\mu},
    \]
    thus we can take
    \[
    L = M + \max_e \fbar_e \left(\frac{2 R + \chi}{\tbar_e \rho}\right)^\mu.
    \]

\end{document}